\documentclass{svmult}

\usepackage{makeidx}         
\usepackage{graphicx}        
\usepackage{multicol}        
\usepackage[bottom]{footmisc}

\makeindex             

\newcommand{\be}{\begin{equation}}
\newcommand{\ee}{\end{equation}}
\newcommand{\ba}{\begin{array}}
\newcommand{\ea}{\end{array}}

\newcommand{\pref}[1]{(\ref{#1})}

\begin{document}

\title*{A comparison between relaxation and Kurganov-Tadmor schemes}
\titlerunning{}
\author{Fausto Cavalli\inst{1} \and Giovanni Naldi\inst{1}  \and
  Gabriella Puppo\inst{2}  \and Matteo   Semplice\inst{1}}
\authorrunning{F. Cavalli, G. Naldi, G. Puppo, M. Semplice}
\institute{Dipartimento di Matematica, Universit\`a di Milano, Via
  Saldini 50, 20133 Milano, ITALY.
  \texttt{\{cavalli,naldi,semplice\}@mat.unimi.it} 
  \and 
  Dipartimento di Matematica, Politecnico di Torino, Corso Duca degli Abruzzi, 24, 10129 Torino, Italy
  \texttt{gabriella.puppo@polito.it}
}
%
%

\maketitle

\begin{abstract}
In this work we compare two semidiscrete schemes for the solution of
hyperbolic conservation laws, namely the relaxation \cite{JX95} and
the Kurganov Tadmor central scheme \cite{KT00}. We are particularly
interested in their behavior under small time steps, in view of future
applications to convection diffusion problems.  The schemes are tested
on two benchmark problems, with one space variable.
\end{abstract}

\section{Motivation}
We are interested in the solution of systems of equations of the form:
\be 
u_t + f_x(u) = D \, p_{xx}(u), \label{eq:general_system}
\ee
where $f(u)$ is hyperbolic, i.e. the Jacobian of $f$ is provided with
real eigenvalues and a basis of eigenvectors for each $u$, while
$p(u)$ is a non decreasing 
Lipshitz continuous function, with Lipshitz constant $\mu$ and $D\geq0$. 

We continue the study of convection diffusion equations with the aid
of high order relaxation schemes started in \cite{arxiv0604572} for
the case of the purely parabolic problem. 

In many
applications, such as multiphase flows in porous media, $p(u)$ is non
linear and possibly degenerate. In these conditions, an implicit
solution of the diffusion term can be computationally very expensive:
in fact it may be necessary to solve large non linear algebraic
systems of equations which, moreover, can be singular at degenerate
points, i.e. where $p(u)=0$. For this reason, it is of interest to
consider the {\em explicit} solution of \pref{eq:general_system}. This
in turn poses one more difficulty. An explicit solution of
\pref{eq:general_system} requires a parabolic CFL condition, that is
stability will restrict the possible choice of the time step $\Delta
t$ to $\Delta t \leq C (\Delta x)^2$, where $\Delta x$ is the grid
spacing. In other words, it may be necessary to choose very small time
steps. But conventional solvers for convective operators typically
work at their best for time steps close to a {\em convective} CFL,
i.e. $\Delta t \leq C \Delta x$. When the time step is much smaller,
they exhibit a very large artificial diffusion of the form $O((\Delta
x)^{2r}/\Delta t)$, where $r$ is the accuracy of the scheme, see for
instance \cite{KT00}. Clearly in these conditions artificial diffusion
becomes very large for $\Delta t \rightarrow 0$. 

As a first step to the numerical solution of problem
\pref{eq:general_system}, we concentrate on semidiscrete schemes for
the solution of the convective part of \pref{eq:general_system}. Such
schemes enjoy an artificial diffusion which depends weakly on $\Delta
t$, and are therefore particularly suited for the solution of
convection-diffusion equations.

We will compare two semidiscrete methods for the
integration of systems of hyperbolic equations. We are interested in
the representation of  solutions which can be characterized by strong
gradients, and, in the degenerate case, even by
discontinuities. Moreover, we are interested in comparing the behavior
of the schemes for small values of $\Delta t$, and for such small
values of the time step, we will investigate the resolution of
discontinuous solutions and the behavior of the error in a few test problems. 

The schemes analyzed in this work are the Kurganov Tadmor central
scheme proposed in \cite{KT00}, and the relaxation scheme proposed in
\cite{JX95}. These methods discretize the equations starting from very
different ideas, however they share some interesting properties. First
of all, they are both semidiscrete schemes. Therefore they require 
separate discretizations in space and time, which is the key to the
fact that artificial diffusion depends mainly on space
discretization. Secondly, they are both Riemann solver free
methods. The Kurganov-Tadmor scheme is based on a central approach:
the solution of the Riemann problem is computed on a staggered cell,
before being averaged back on the standard grid. In this fashion, the
numerical solution is updated on the edges of the staggered grid,
where it is smooth, and can be computed via a Taylor expansion, with
no need to solve the actual Riemann problem. The relaxation scheme
instead moves the non linearities of the convective term to a stiff
source term, and the transport part of the system becomes linear, with
a fixed and well known characteristic structure. Thus again there is
no need to use approximate or exact Riemann solvers.  

For these reasons both schemes can be applied as black-box
methods to a fairly general class of balance laws.

\section{Results}

For the Kurganov-Tadmor (KT) scheme we have followed the componentwise
implementation of the method described in \cite{KT00}. The scheme is
written in conservation form, with numerical flux: 
\be
\ba{ll}
F_{j+1/2}(t) =  \frac12 & \left[ f(u^+_{j+1/2}(t)) + f(u^-_{j+1/2}(t)) \right.\\
 & \left. -a_{j+1/2}(t) \left( u^+_{j+1/2}(t) - u^-_{j+1/2}(t) \right) \right], 
\ea
\ee
where $u^+_{j+1/2}(t)$ and $u^-_{j+1/2}(t) $ are the boundary
extrapolated data, computed at the edges of each cell with a piecewise
linear reconstruction at time $t$, and $a_{j+1/2}(t)$ is a measure of
the maximum propagation speed at the cell edge. This value for the
case of systems of equations, in particular in the non convex case,
must be carefully tuned, and it is the same for all components, when
the scheme is implemented componentwise. 

On the other hand, the relaxation scheme  requires an accurate choice of the
subcharacteristic velocities $A^2$. The relaxation system is
\begin{equation} \label{eq:rel}
\left\{
\begin{array}{ll}
\displaystyle
\frac{\partial \vec{u}}{\partial t} + \frac{\partial \vec{v}}{\partial x} = 0 \\
\\
\displaystyle 
\frac{\partial \vec{v}}{\partial t} + A^2 \frac{\partial \vec{u}}{\partial x}  = 
    -\frac1\varepsilon \left(\vec{v}-f(\vec{u})\right) 
\end{array} 
\right.   .
\end{equation}
For a scalar conservation law, we take $A^2=\max(|f'(u)|)$ as in
\cite{JX95}, while for the Euler system of gas-dynamics
we take $A^2$ to be the diagonal matrix with entries
$\max_j(|u_j-c_j|),\max_j(|u_j|),\max_j(|u_j+c_j|)$.
Here $u$ is the velocity and $c$ the speed of sound. We update these
quantities at each time step, so that $A^2$ can be chosen as small as
possible (in the paper \cite{JX95} $A^2$ was chosen as a constant
diagonal matrix but this results in a larger numerical diffusion). 

Due to the diagonal form of $A^2$, the convective operator is block
diagonal with $2\times2$ blocks. Each block is independently
diagonalized and we compute the numerical fluxes using a second order
ENO reconstruction \cite{HEOC87}.

As $\varepsilon\rightarrow 0$, the system \pref{eq:rel} formally
relaxes to the original conservation laws, provided the
subcharacteristic condition holds, namely that $ (A^2 - (f'(u))^2)$ is
positive-definite.

We use the second order Heun Runge-Kutta method for the time
integration of both the KT and the relaxation schemes.

\begin{table}
\begin{center}
\begin{tabular}{|r|c|c|c|c|}
\hline
&\multicolumn{2}{c|}{\ Convective CFL\ }&\multicolumn{2}{c|}{\ Parabolic CFL\ }\\
\hline
&KT&Relax&KT&Relax\\
\hline
$\;20\;$   &  2.03E-1      & 2.16E-1 & $\;$6.19E-1$\;$ & $\;$1.02E-1$\;$\\
$\;40\;$   &  7.58E-2      & 7.66E-2 & 2.04E-1 & 4.58E-2\\
$\;80\;$   &  2.71E-2      & 2.73E-2 & 9.10E-2 & 1.34E-2\\
$\;160\;$  & 8.22E-3      & 8.25E-3 & 2.67E-2 & 3.82E-3\\
$\;320\;$  & 2.29E-3  & 2.29E-3 & 7.62E-3 & 1.03E-3\\
$\;640\;$  & 6.11E-4  & 6.12E-4 & 2.06E-3 & 2.77E-4\\
$\;1280\;$ & $\;$1.61\mbox{E}-4$\;$ & $\;$1.61E-4$\;$ & & \\
\hline
\end{tabular}
\end{center}
\caption{Linear advection of a sine function. Errors in $\mathrm{L}^1$
  at $t=1$.}
\label{tab:conv}
\end{table}

Table \ref{tab:conv} shows the errors in the $\mathrm{L}^1$ norm for
the linear advection equation $u_t+u_x=0$ with initial data 
$u(x,0)=\sin(2\pi{}x)$. We use the standard convective CFL condition
$\Delta t= C \Delta x$ and the parabolic CFL, $\Delta t= C (\Delta
x)^2$. We note that the errors are almost the same for the two schemes
for the convective CFL, while the relaxation scheme seems superior for
the parabolic CFL.

\begin{figure}
\begin{center}
\includegraphics[height=.35\textwidth,width=.45\textwidth]{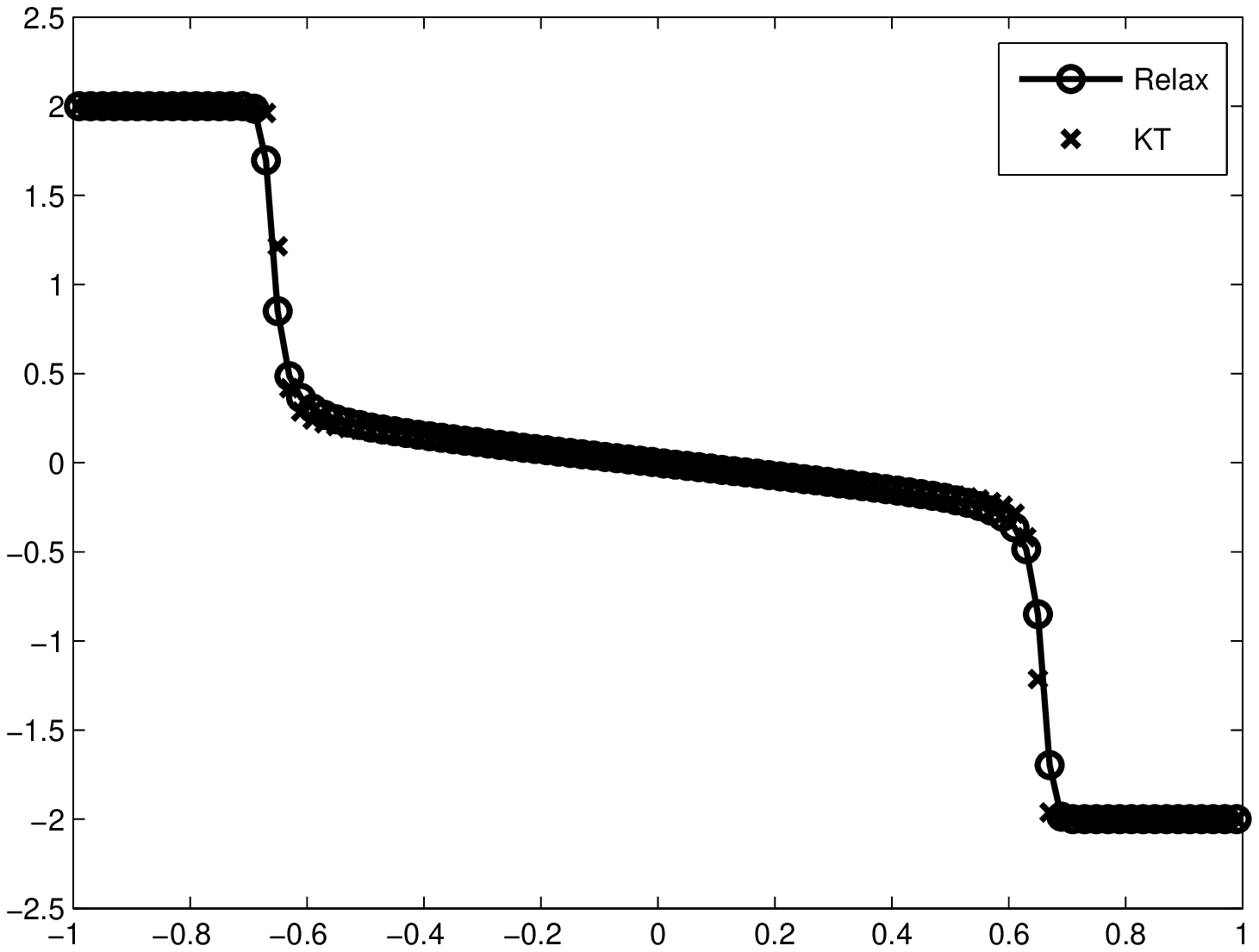}
\hfil
\includegraphics[height=.35\textwidth,width=.45\textwidth]{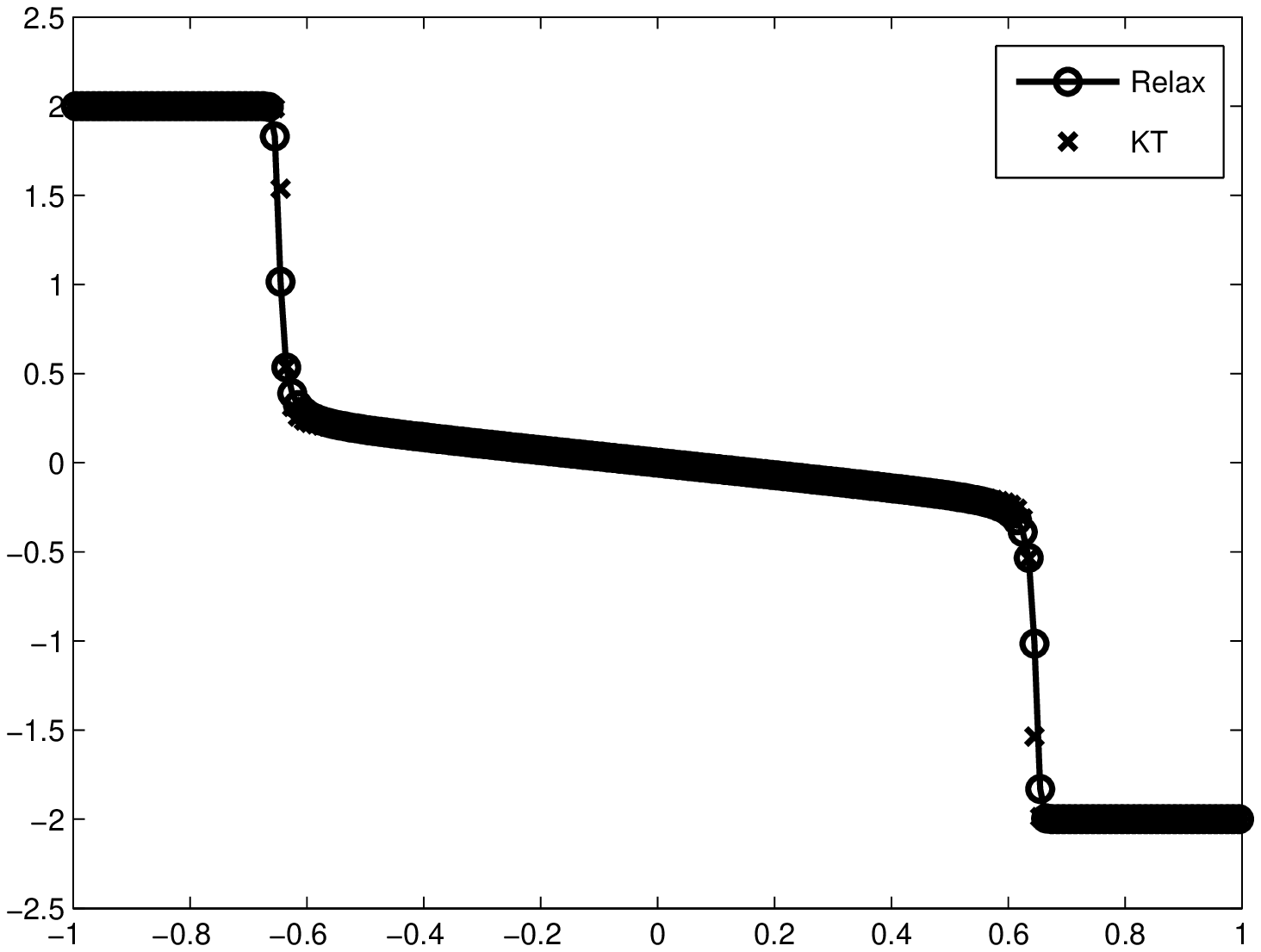}
\end{center}
\caption{Non-convex flux. Kurganov-Tadmor and relaxation schemes, with
$n=100$ (left) and $n=200$ (right).}
\label{fig:nonconv}
\end{figure}

A key requirement for a numerical scheme for conservation laws is the
ability to pick the entropy solution in non-convex problems. Here we
show a Riemann problem for the non-convex flux
$f(u)=(u^2-1)(u^2-4)/4$, as in \cite{KT00}. The Riemann problem breaks
into two shocks connected by a rarefaction wave. The results are shown
in Figure \ref{fig:nonconv}. Clearly both schemes are able to resolve
the correct discontinuities and they have approximately the same
resolution, the KT scheme being slightly less diffusive.

\begin{figure}
\begin{center}
\includegraphics[height=.45\textwidth,width=.45\textwidth]{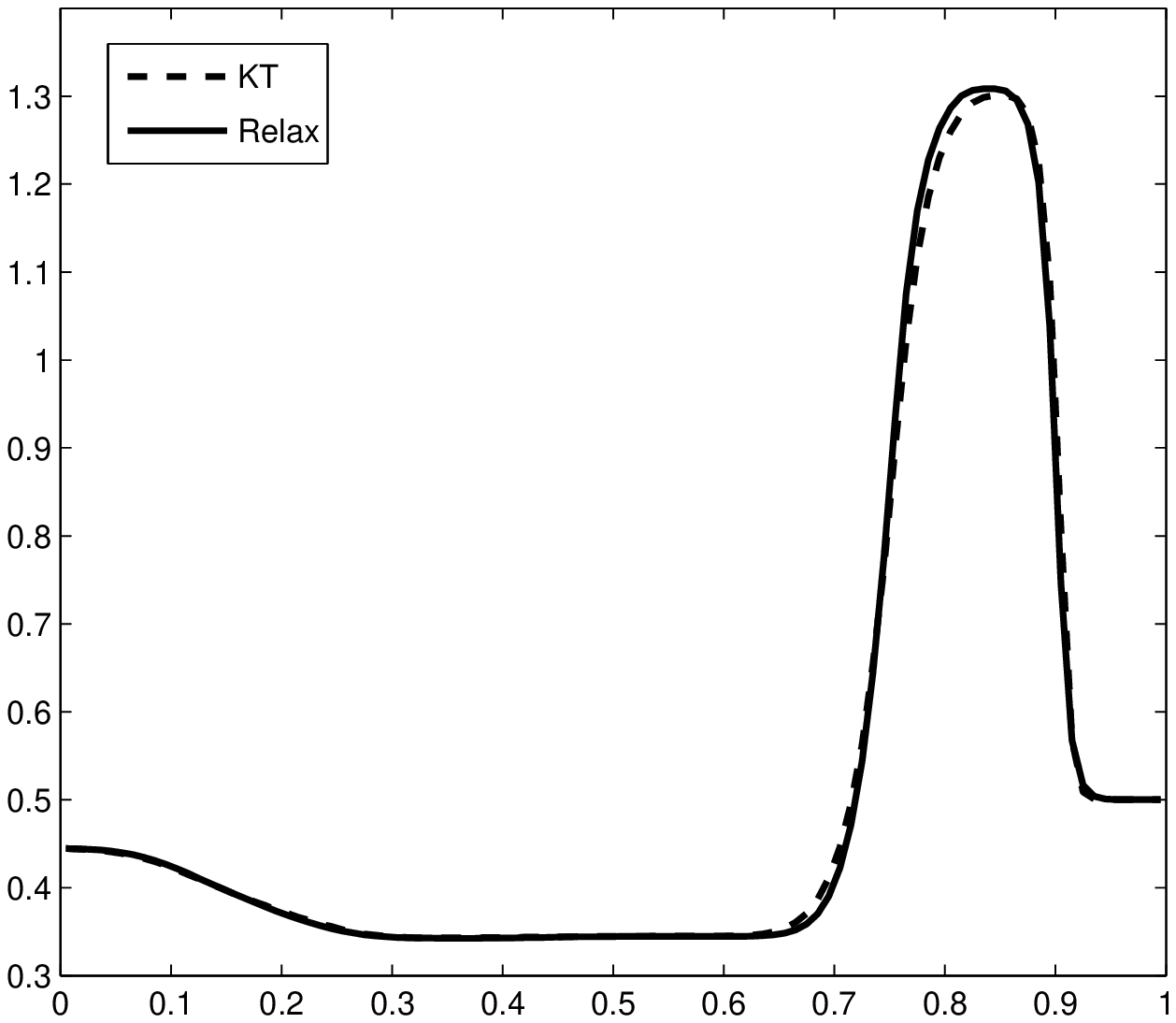}
\hfil
\includegraphics[height=.45\textwidth,width=.45\textwidth]{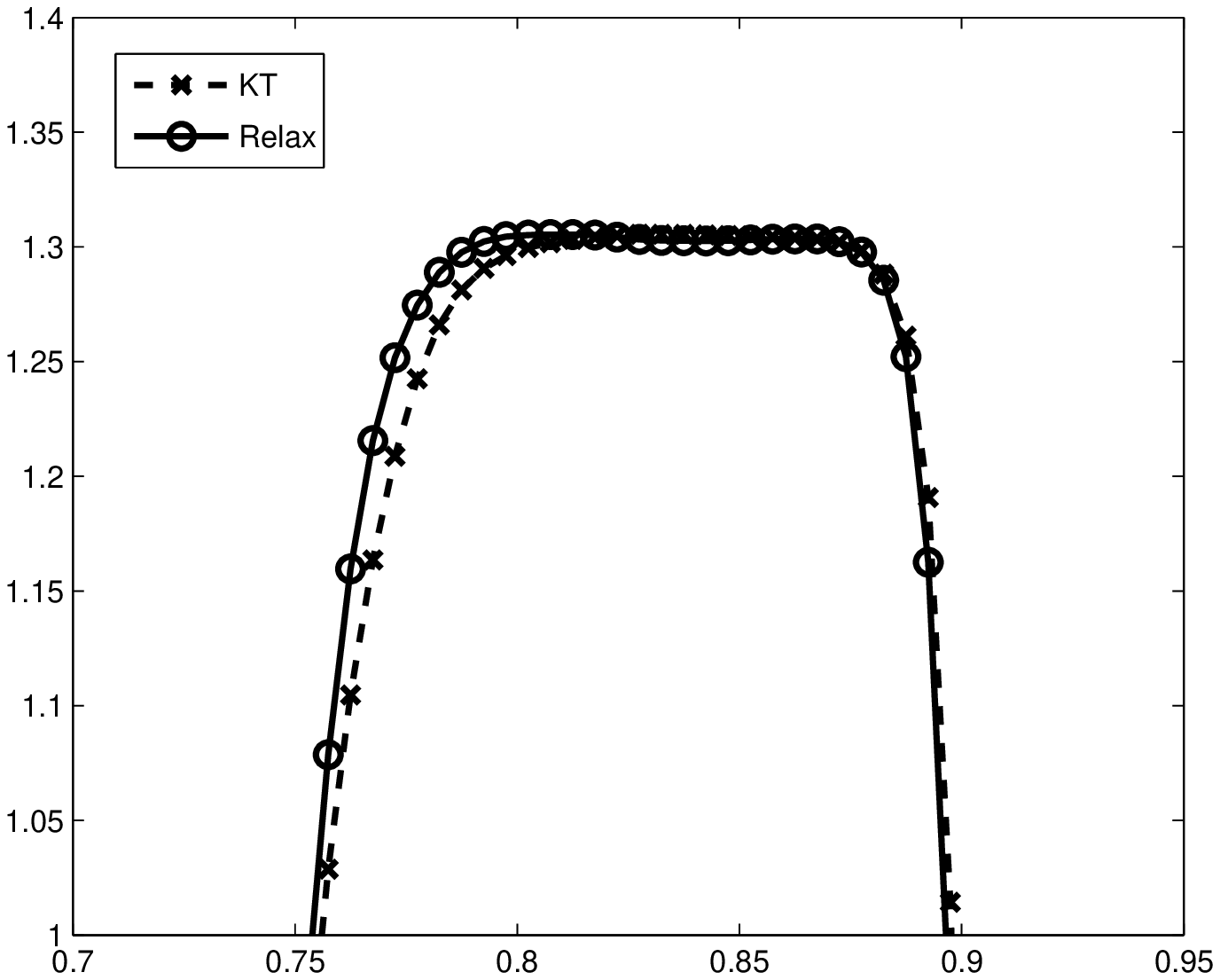}
\end{center}
\caption{Lax Riemann problem, density component. KT (dashed) and
  relaxation (solid line) with $n=100$ (left) and $n=200$ (right),
  where a detail of the density peak is shown.}
\label{fig:lax}
\end{figure}

Figure \ref{fig:lax} shows the density component of the Lax Riemann
problem in gas dynamics. The accurate choice suggested above for
the matrix $A^2$ in the relaxation system yields a slightly higher
resolution than KT.

\noindent{\bf Concluding remarks}

\noindent
We have compared two semidiscrete schemes for conservation laws. We
find that although the schemes are constructed with very different
philosophies, they yield comparable results on some significant test
problems. We think that the relaxation scheme is slightly more robust,
since it results from the relaxation of a viscous profile, provided
the subcharacteristic condition is satisfied. Also, the actual errors
obtained with a parabolic CFL in Table \ref{tab:conv} seem to favour
the relaxation scheme.

We also wish to mention higher order extensions of the schemes studied
in this work: namely the third order central upwind scheme described
in \cite{KNP01}, endowed with a more carefully crafted artificial
diffusion with respect to \cite{KT00} and the third order extension of
the relaxation scheme proposed in \cite{Sea06}.


%
%
 \bibliographystyle{alpha}
 \bibliography{ECMI06}
%


\printindex
\end{document}